\documentclass[]{interact}
\usepackage{subfigure}   

\usepackage[numbers,sort&compress]{natbib}
\usepackage{enumerate}
\usepackage{graphicx}
\usepackage{amsmath} 
\usepackage{amssymb}
\usepackage{algorithm,setspace}
\usepackage{algpseudocode}
\DeclareMathOperator*{\argmin}{argmin}
\DeclareMathOperator{\inte}{int}
\DeclareMathOperator{\dom}{dom}
\DeclareMathOperator{\rng}{range}
\DeclareMathOperator{\prox}{{prox}}
\bibpunct[, ]{[}{]}{,}{n}{,}{,}
\makeatletter
\def\NAT@def@citea{\def\@citea{\NAT@separator}}
\makeatother

\theoremstyle{plain}
\newtheorem{theorem}{Theorem}[section]
\newtheorem{lemma}[theorem]{Lemma}

\newtheorem{proposition}[theorem]{Proposition}

\theoremstyle{definition}

\theoremstyle{remark}
\newtheorem{remark}{Remark}

\begin{document}


\title{Modified Bregman Golden Ratio Algorithm for Mixed Variational Inequality Problems}

\author{
\name{Gourav Kumar\textsuperscript{a}\thanks{Gourav Kumar,  Email: ma24r001@smail.iitm.ac.in; V.Vetrivel, Email:vetri@iitm.ac.in} and V. Vetrivel\textsuperscript{a}}
 \affil{\textsuperscript{a}Department of Mathematics, Indian Institute of Technology Madras, Chennai,\\ ~~600036, India}}

\maketitle

\begin{abstract}
In this article, we provide a modification to the Bregman Golden Ratio Algorithm (B-GRAAL). We analyze the B-GRAAL algorithm with a new step size rule, where the step size increases after a certain number of iterations and does not require prior knowledge of the global Lipschitz constant of the cost operator. Under suitable assumptions, we establish the global iterate convergence as well as the R-linear rate of convergence of the modified algorithm. 
The numerical performance of the proposed approach is validated for the matrix game problem and the sparse logistic regression problem in machine learning.  
\end{abstract}
\par

\section{Introduction}
Let $\mathcal{H}$ be a real finite dimentional Hilbert space. A mixed variational inequality problem (MVIP) is a problem to find $w^*\in\mathcal{H}$ such that 
\begin{equation}\label{eqn_defvi}
    \langle A(w^*),w-w^*\rangle+g(w)-g(w^*)\geq 0~\forall w\in \mathcal{H}, 
\end{equation}
where $A:\mathcal{H}\to\mathcal{H}$ is the cost operator, and $g:\mathcal{H}\to\mathbb{R}\cup\{+\infty\}$ is an extended real-valued function. We denote the set of solutions of (\ref{eqn_defvi}) by $S$. Mixed variational inequality problems (MVIPs) of the form (\ref{eqn_defvi}) generalize various optimization problems encountered in nonlinear programming and variational analysis, including minimization problems, linear complementarity problems, and variational inequalities, see \cite{addi2017complementarity,Ansari,fukushima1992equivalent,popov}. These problems have widespread applications in diverse fields, such as data science, image processing, mechanics, control theory, economics, structural engineering, and many more; see \cite{addi2017complementarity,daniele2003equilibrium,grad2021solving,han1995finite,jolaoso2022inertial,ju2022solving,ju2021novel,ju2024new}, and the references therein.\par
Many algorithms solving the variational inequality (\ref{eqn_defvi}) (for example; see \cite{fukushima1992equivalent,marcotte1995convergence,pang1982iterative}), need the exact value of the global Lipschitz constant of the cost operator $A$. This assumption is often overly restrictive, as determining the Lipschitz constant of $A$ can be more challenging than solving the original problem itself. Furthermore, even when the Lipschitz constant $L$ is known, using a constant step size rule that inversely depends on $L$ can be highly restrictive. In particular, when $L$ is large, the step size $ \lambda$ is restricted to a very small value, which greatly limits the possible choices for $\lambda$ and may slow down convergence for a given number of iterations. To overcome these issues, several researchers proposed step size sequences which do not require the value of Lipschitz constant, rather these algorithms approximate the local Lipschitz constant at each iterate. To generate such step size sequences the authors used \emph{backtracking} procedures (see \cite{bello2015variant,censor1998interior,iusem1997variant,malitsky2018proximal,malitsky2020forward}).
Although these methods address the aforementioned drawbacks, they can become costly in terms of overall run time due to the potentially large number of steps required for the backtracking procedure in each iteration. \par
To address this issue, Malitsky \cite{malitsky2020golden} proposed an adaptive strategy for step size selection. Subsequently, several authors adopted the adaptive step size technique to address the limitations of the backtracking procedure \cite{alacaoglu2020convergence,soe,zhou2024generalized}. To build upon earlier methods, Hoai \cite{hoai2024new} recently proposed an algorithm to solve (\ref{eqn_defvi}), in which the step size sequence increases after a finite number of iterations. One advantage of using an increasing step size is that it allows for adaptation when the initial step size is too small, particularly in cases where the Lipschitz constant is very large. If the step size begins at a low value, an increasing sequence can gradually adjust it over iterations, improving efficiency. This flexibility provides a key advantage over fixed or nonincreasing step size strategies.\par
Another approach to improve algorithms to solve (\ref{eqn_defvi}) is to replace the Euclidean distance  by other nonlinear distances such as Bregman distance. The usage of Bregman distance and Bregman proximal operators provides flexibility when deciding which projection to compute. For example, in convex-concave matrix games, the Bregman proximal operator appears as a projection onto the probability simplex, that has a straightforward closed form expression with regard to the \emph{Kullback-Leibler (KL)} divergence but not with respect to the conventional Euclidean distance. For further advantages of the Bregman distance over the Euclidean distance, we refer the reader to \cite{chen2011projection,dai2022distributed,wang2013projection,tam2023bregman}. The extant literature has several approaches for solving (\ref{eqn_defvi}) with Bregman proximal operators \cite{censor1998interior,gibali2018new,gibali2020fast,van2022modified,hieu2022two,jolaoso2022single,jolaoso2020weak,jolaoso2020two,nomirovskii2019convergence,semenov2021adaptive}.\par
Most of these methods assume a Lipschitz condition without requiring the explicit value of the Lipschitz constant. Nevertheless, they often rely on a backtracking procedure. 
Recently, Tam and Uteda \citep[Algorithm 2]{tam2023bregman}, studied the Bregman modifications of the \emph{Golden Ratio Algorithm}, which iterates as follows. 
\begin{algorithm}[H]
\caption{The Bregman-Golden Ratio Algorithm (B-GRAAL)}\label{algo_01} 
\begin{algorithmic}[1] 
\Statex \textbf{Input}: $w_1,\bar{w}_0\in\inte\dom h$ and a step size $\lambda\in\left(0,\frac{\sigma\phi}{2L}\right]$.
\Statex\textbf{for} k=1,2,\ldots,
\Statex Compute the next iterations:
\begin{align*}
    \bar{w}_k=&~(\nabla h)^{-1}\left(\frac{(\phi-1)\nabla h(w_k)+\nabla h(\bar{w}_{k-1})}{\phi}\right),\\w_{k+1}=&~\argmin_{w\in \mathcal{H}}\left\{\langle A(w_k),w-w_k\rangle + g(w)+\frac{1}{\lambda_k}B_{h}(w,\bar{w}_k)\right\}.
\end{align*}
\end{algorithmic}
\end{algorithm}\noindent
This method does not require the backtracking procedure. However, to implement the algorithm, one needs to know the exact value of the global Lipschitz constant. To overcome this issue, in the same article \citep[Algorithm 3]{tam2023bregman}, the authors proposed the \emph{Bregman adaptive Golden Ratio Algorithm (B-aGRAAL)}. Implementation of the B-aGRAAL algorithm requires neither backtracing nor the value of the Lipschitz constant.\par 
Motivated by the work of Tam and Uteda \cite{tam2023bregman}, in this article, we propose a modification to the B-GRAAL algorithm. We study the B-GRAAL algorithm with a new step size rule. Our choice of step size sequence is motivated by the step size used by Hoai in \cite{hoai2024new}. This choice of step size sequence eliminates the need to know the value of global Lipschitz constant to implement the B-GRAAL algorithm. Furthermore, unlike previous approaches that utilize Bregman divergence to solve MVIP (\ref{eqn_defvi}), our step size sequence increases after a finite number of iterations.
We present a comprehensive convergence analysis of the proposed algorithm. Further, we establish the R linear rate of convergence of the algorithm. The numerical results given in Section \ref{sec:5} reveal that the proposed algorithm significantly improves the performance of the B-GRAAL algorithm.\par
This paper is structured as follows. Section \ref{sec:2} begins with the essential definitions and fundamental results as a foundation for the convergence analysis. Section \ref{sec:3} introduces the main results, including the proposed algorithm and its convergence analysis. Section \ref{sec:5} presents numerical results to illustrate the effectiveness of the proposed method.
\section{Preliminaries}\label{sec:2}
\noindent
In this article, $\mathcal{H}$ represents a real, finite-dimensional Hilbert space equipped with an inner product $\langle \cdot, \cdot\rangle$ and the corresponding norm $\lVert\cdot\rVert$.
Given $g:\mathcal{H}\rightarrow (-\infty, +\infty]$, domain of $g$ is represented by $\dom g:=\{x\in\mathcal{H}:g(x)<+\infty\}.$ The subdifferential of $g$ at $x\in\dom g$ is given by
\[\partial g(x):=\{w\in\mathcal{H}:g(x)-g(y)-\langle w, x-y\rangle\leq 0~\forall y\in\mathcal{H}\}.\]
\begin{lemma}\emph{\citep[Lemma 2.3]{hoai2024new}}\label{lm_pre_1}
    Let $(x_k)$ and $(y_k)$ be two sequences of nonnegative numbers satisfying:
    \[x_{k+1}\leq x_k-y_k~\forall k\in\mathbb{N}.\]
    Then, $(x_k)$ is convergent and $\sum_{k=0}^{+\infty} y_k<+\infty.$
\end{lemma}
For $h:\mathcal{H}\rightarrow (-\infty, +\infty]$, the Bregman distance with respect to $h$ is the function $B_{h}:\mathcal{H}\times \inte\dom h\rightarrow (-\infty,+\infty]$ defined as 
\[B_{h}(x,y):=h(x)-h(y)-\langle\nabla h(y), x-y\rangle.\]
For $\alpha$-strongly convex function $h$, $B_{h}$ satisfies
\[B_{h}(x,y)\geq\frac{\alpha}{2}\lVert x-y\rVert^{2}~\forall (x,y)\in \mathcal{H}\times \inte\dom h.\]
The Bregman distance with respect to two different choices of $h$ is given as follows \cite{bauschke2009bregman,van2022modified}. Let $z=(z_1,z_2,\ldots,z_n)^\top$ and $w=(w_1,w_2,\ldots,w_n)^\top$ be two points in ${\mathbb{R}}^n.$ 
\begin{enumerate}[(i)]
    \item The Kullback-Leibler distance is given as
    \[B^{KL}_h(z,w):=\sum\limits_{{i=1}}^n z_i\left(\log\left(\frac{z_i}{w_i}\right)-1\right)+\sum\limits_{i=1}^n w_i,\]
    which is induced by the function $h(z)=\sum\limits_{i=1}^n z_i\log z_i$,  with the domain $\dom h^{KL}=\{z\in{\mathbb{R}}^n:z_i>0\}.$
    \item The squared Mahalanobis  distance is given as 
    \[B^{SM}_h(z,w)=\frac{1}{2}(z-w)^\top Q(z-w),\]
    which is induced by the function $h^{SM}(z)=\frac{1}{2}z^\top Q z$, with $\dom h^{SM}={\mathbb{R}}^n$, where $Q=\text{diag}(1,2,\ldots,n).$
\end{enumerate}
\begin{remark}
    If we take $h=\frac{1}{2}\lVert .\rVert^2$. Then, the Bregman distance simplifies to 
    $B_h(z,w)=\frac{1}{2}\lVert z-w\rVert^2$, 
    which is the squared Euclidean distance.
\end{remark}
In this work, we use the following properties of the Bregman distance.
\begin{proposition}\emph{\citep[Proposition 2.1]{tam2023bregman}}\label{prp_1} For the Legendre function $h:\mathcal{H}\rightarrow (-\infty, +\infty]$, the statements listed below are true.
\begin{enumerate}[(i)]
    \item\label{prp_1_itm_1} For all $x,y\in\inte\dom h$ and $w\in\dom h$,
    \[B_{h}(w,x)-B_{h}(w,y)-B_{h}(y,x)=\langle\nabla h(x)-\nabla h(y),y-w\rangle.\]
    \item\label{prp_1_itm_2} Let $x\in\dom h, y,z, w\in\inte\dom h$ with $\beta\in\mathbb{R}$ and $\nabla h(y)=\beta\nabla h(z)\\+(1-\beta)\nabla h(w)$. Then, 
    \[B_{h}(x,y)=\beta\Big(B_{h}(x,z)-B_{h}(y,z)\Big)+(1-\beta)\Big(B_{h}(x,w)-B_{h}(y,w)\Big).\]
\end{enumerate}
\end{proposition}\noindent
The (left) Bregman proximal operator of a function $f:\mathcal{H}\rightarrow (-\infty,\infty]$ is the operator given by
\[{\prox}^h_f(y):=\argmin_{x\in\mathcal{H}}\{f(x)+B_{h}(x,y)\}~\forall y\in\inte\dom h.\]
Since this work exclusively requires the left Bregman proximal operator, we will henceforth refer to it simply as the Bregman proximal operator.
Next, we state some properties of the Bregman proximal operator useful in our algorithm analysis.
\begin{proposition}\emph{\citep[Proposition 2.2]{tam2023bregman}}\label{prp_2}
    Let $f:\mathcal{H}\rightarrow (-\infty, +\infty]$ be proper, convex, lower semicontinuous, and let $h:\mathcal{H}\rightarrow (-\infty, +\infty]$ be Legendre with $\inte\dom h\cap \dom f\neq\emptyset$.
    \begin{enumerate}[(i)]
        \item $\rng(\prox^{h}_{f})\subseteq\inte\dom h\cap\dom f.$
        \item $\prox^{h}_{f}$ is single-valued on $\dom(\prox^{h}_{f})\subseteq\inte\dom h$. In addition, if $h+f$ is super-coercive, then $\dom(\prox^{h}_{f})=\inte\dom h.$
        \item\label{prp2_itm_3} Let $y\in\dom(\prox^{h}_{f})$. Then, $x=\prox^{h}_{f}(y)$ if and only if, for all $z\in \mathcal{H}$, we have
        \[f(x)-f(z)\leq\langle\nabla h(y)-\nabla h(x), x-z \rangle.\]
    \end{enumerate}
\end{proposition}
\section{The Modified Bregman-Golden Ratio Algorithm}\label{sec:3}\noindent
In this section, we present the modified Bregman Golden Ratio algorithm and discuss its convergence. We take the following assumptions for our subsequent convergence analysis.
\allowdisplaybreaks
\begin{enumerate}[{A}.1]
 \item\label{assum_1} $g:\mathcal{H}\rightarrow (-\infty, +\infty]$ is proper, convex, and lower semicontinuous.
    \item\label{assum_2} $h:\mathcal{H}\rightarrow (-\infty, +\infty]$ is continuously differentiable, Legendre, and $\alpha$-strongly convex. Further, we assume that for any sequence $(w_n)\subseteq \inte\dom h$ converging to some $w\in \dom h,$  $B_{h}(w, w_n)\rightarrow 0$.
    \item\label{assum_3} $A:\mathcal{H}\rightarrow \mathcal{H}$ is monotone and $L$-Lipschitz continuous on $\dom g \cap \inte\dom h\neq\emptyset$. 
    \item\label{assum_4} $\mathcal{S}:=S\cap\dom h\neq\emptyset$.
\end{enumerate}
Let $\phi:= \frac{\sqrt{5}+1}{2}$ represent the Golden Ratio that follows the identity ${\phi}^2=\phi +1.$ The proposed algorithm is formally outlined in Algorithm \ref{algo_1}.
\vspace{-10pt}
\begin{algorithm}[H]
\caption{The Modified B-GRAAL for (\ref{eqn_defvi})}  \label{algo_1} 
\begin{algorithmic}[1] 
\Statex \textbf{Initialization}: Let $\lambda_0>0$, $0<\eta_1<\eta_0<\frac{\phi}{2}$ and a sequence $\gamma_k$ of positive real numbers satisfying
\begin{equation*}
    \sum_{k=0}^{+\infty}\gamma_k< +\infty.
\end{equation*}
Choose $w_0=\bar{w}_0\in \mathcal{H}$, $w_1\in \mathcal{H}$ and set $k=1$.
\Statex\textbf{Step 1}: If 
\begin{equation*}
    \lVert A(w_k)-A(w_{k-1})\rVert>\frac{\eta_0\alpha}{\lambda_{k-1}}\lVert w_k-w_{k-1}\rVert,
\end{equation*}
then
\begin{equation}\label{eqn_roc_2}
    \lambda_k=\eta_1\alpha\frac{\lVert w_k-w_{k-1}\rVert}{\lVert A(w_k)-A(w_{k-1})\rVert}
\end{equation}
else 
\begin{equation}\label{eqn_roc_3}
    \lambda_k=(1+\gamma_{k-1})\lambda_{k-1}.
\end{equation}
\Statex \textbf{Step 2}: Compute
\begin{align}
    \bar{w}_k=&~(\nabla h)^{-1}\left(\frac{(\phi-1)\nabla h(w_k)+\nabla h(\bar{w}_{k-1})}{\phi}\right)\label{eqn_stp2_1},\\w_{k+1}=&~\argmin_{w\in \mathcal{H}}\left\{\langle A(w_k),w-w_k\rangle + g(w)+\frac{1}{\lambda_k}B_{h}(w,\bar{w}_k)\right\}\label{eqn_step_2}.
\end{align}
Set $k:=k+1$ and return to \textbf{Step 1}.
\end{algorithmic}
\end{algorithm}\noindent
\noindent
It is interesting to note that (\ref{eqn_step_2}) in Algorithm \ref{algo_1}, can be equivalently written in the following forms \citep[Remark 3.2]{tam2023bregman}: 
\begin{align*}
    w_{k+1} =& \prox^{h}_{\lambda_k f}(\bar{w}_k), ~\text{where}~f(w) = \langle A(w_k),w - w_k \rangle + g(w)\\=&\prox^{h}_{\lambda_k g} \left( (\nabla h)^{-1} (\nabla h(\bar{w}_k) - \lambda_k A(w_k)) \right).
\end{align*}
By applying a proof analogous to that of Lemma 3.1 in \cite{hoai2024new}, it can be established that the step size sequence $(\lambda_k)$ satisfies the following properties.
\begin{lemma}\label{mm_lm}
Suppose $A$ is $L$-Lipschitz continuous and $(\lambda_k)$ is the step size sequence produced by {Algorithm} \ref{algo_1}. Then, 
    \begin{enumerate}[(i)]
        \item\label{lem_stpsize_prp1} $\lambda_k\geq\min\left\{\frac{\eta_1}{L},\lambda_0\right\}>0~\forall$ $k\in\mathbb{N}$;
        \item $(\lambda_k)$ is convergent;
        \item\label{lm_prt_3} there exists $k_1\in\mathbb{N}$ satisfying $\lambda_{k+1}\geq\lambda_k$ for all $k\geq k_1.$ Equivalently, there exists $k_1\in\mathbb{N}$ such that
        \[\lVert A(w_k)-A(w_{k-1})\rVert\leq\frac{\eta_0\alpha}{\lambda_{k-1}}\lVert w_k-w_{k-1}\rVert~\forall k\geq k_1.\]
    \end{enumerate}
\end{lemma}
\begin{remark}\label{obs_1}
    Since the positive series $\sum_{k=0}^{+\infty}\gamma_k< +\infty$ in Algorithm \ref{algo_1} converges and $0<\eta_1<\eta_0<\frac{\phi}{2}$, there exists $\bar{k}\in\mathbb{N}$ such that 
\[\gamma_{k-1}<\frac{\phi}{\eta_0}-1~\forall k\geq \bar{k}.\]
\end{remark}
Next, we present two key lemmas that play a fundamental role in proving the convergence of the algorithm.
\begin{lemma}\label{lm_mm2}Suppose that Assumption A.\ref{assum_2} holds. If $(w_k)$ is generated by Algorithm \ref{algo_1}, then there exists $\bar{k}\in\mathbb{N}$ such that
    \[\lambda_k\langle A(w_k)-A(w_{k-1}), w_k-w_{k+1}\rangle<\frac{\phi}{2}\Big(B_{h}(w_k,w_{k-1})+B_{h}(w_{k+1},w_k)\Big)~\forall k\geq \bar{k}.\]
\end{lemma}
\allowdisplaybreaks
{\it \textbf{Proof}}
    We consider the following two possible cases.
    \begin{enumerate}[{Case} 1.]
        \item\label{lm_cs_1}  $\lVert A(w_k)-A(w_{k-1})\rVert>\frac{\eta_0\alpha}{\lambda_{k-1}}\lVert w_k-w_{k-1}\rVert.$\\
        By using Cauchy-Schwarz inequality, relation (\ref{eqn_roc_2}), definition of $\lambda_k$ and $\alpha$- strong convexity of $h$, we get
        \allowdisplaybreaks
       \begin{align*}
                &\lambda_k\langle A(w_k)-A(w_{k-1}),w_k-w_{k+1}\rangle\\& \leq \lambda_k\lVert A(w_k)-A(w_{k-1})\rVert~ \lVert w_k-w_{k+1}\rVert\\&=\eta_1 \alpha\lVert w_k-w_{k-1}\rVert ~\lVert w_k-w_{k+1}\rVert\\&<\frac{\phi}{2}\alpha\lVert w_k-w_{k-1}\rVert~ \lVert w_k-w_{k+1}\rVert\\&\leq\frac{\phi}{2}\frac{\alpha}{2}\Big(\lVert w_k-w_{k-1}\rVert^2+\lVert w_k-w_{k+1}\rVert^2\Big)\\ &\leq\frac{\phi}{2}\Big(B_{h}(w_k,w_{k-1})+B_{h}(w_{k+1},w_k)\Big).
        \end{align*}
        \item \label{lm_cs_2} $ \lVert A(w_k)-A(w_{k-1})\rVert\leq\frac{\eta_0\alpha}{\lambda_{k-1}}\lVert w_k-w_{k-1}\rVert$.\\
        By applying the Cauchy-Schwarz inequality, relation (\ref{eqn_roc_3}), Remark \ref{obs_1} and $\alpha$- strong convexity of $h$, we obtain the following $\forall k\geq\bar{k}$,
       \begin{align*}
                &\lambda_k\langle A(w_k)-A(w_{k-1}),w_k-w_{k+1}\rangle\\&\leq \lambda_k\lVert A(w_k)-A(w_{k-1})\rVert~ \lVert w_k-w_{k+1}\rVert\\ &\leq \lambda_k\frac{\eta_0\alpha}{\lambda_{k-1}}\lVert w_k-w_{k-1}\rVert~\lVert w_k-w_{k+1}\rVert\\&=\eta_0(1+\gamma_{k-1})\alpha\lVert w_k-w_{k-1}\rVert~\lVert w_k-w_{k+1}\rVert\\& <\frac{\phi}{2}\frac{\alpha}{2}\Big(\lVert w_k-w_{k-1}\rVert^2+\lVert w_k-w_{k+1}\rVert^2\Big)\\&\leq\frac{\phi}{2}\Big(B_{h}(w_k,w_{k-1})+B_{h}(w_{k+1},w_k)\Big).
        \end{align*}
    \end{enumerate}
     Thus, from Case \ref{lm_cs_1} and Case \ref{lm_cs_2}, we have
\[\lambda_k\langle A(w_k)-A(w_{k-1}), w_k-w_{k+1}\rangle<\frac{\phi}{2}\Big(B_{h}(w_k,w_{k-1})+B_{h}(w_{k+1},w_k)\Big)~\forall k\geq \bar{k}.\]
\qedhere
\allowdisplaybreaks
\begin{lemma}\label{lm_1}
    Suppose that Assumptions A.\ref{assum_1}-A.\ref{assum_4} hold.  
    Let $w^*\in \mathcal{S}$ be arbitrary. Then, there exists $k'\in\mathbb{N}$ such that for all $k\geq k'$ the sequences $(w_k)$ and $(\bar{w}_k)$ generated by {Algorithm \ref{algo_1}} satisfy 
    \allowdisplaybreaks
    \begin{align}\label{eqn_lmprf_main}
            0 &\leq (\phi+1)B_{h}(w^*,\bar{w}_{k+1})+\frac{\phi}{2}B_{h}(w_{k+1},w_k)-\phi B_{h}(w_{k+1},\bar{w}_{k+1})\nonumber\\
            &\leq (\phi+1)B_{h}(w^*,\bar{w}_k)+\frac{\phi}{2}B_{h}(w_k,w_{k-1})-\phi B_{h}(w_k,\bar{w}_k)\nonumber\\&~~~~~-\left(1-\frac{3}{2\phi}\right) B_{h}({w}_{k+1},\bar{w}_k).
    \end{align}
\end{lemma}
\allowdisplaybreaks
{\it \textbf{Proof}}
    From (\ref{eqn_step_2}), we have
    \begin{equation}\label{eqn_lmprf_1}
        w_{k+1}=\argmin_{w\in \mathcal{H}}\left\{\langle A(w_k),w-w_k\rangle + g(w)+\frac{1}{\lambda_k}B_{h}(w,\bar{w}_k)\right\}.
    \end{equation}
 Equation (\ref{eqn_lmprf_1}) can be equivalently written in terms of Bregman proximal operator as follows:
    \[w_{k+1}=\prox^{h}_{\lambda_{k}f}(\bar{w}_k), \text{where}~f(w)=\langle A(w_k),w-w_k\rangle+g(w).\]
    By applying Proposition \ref{prp_2}(\ref{prp2_itm_3}) with $z:=w\in \dom h\cap\dom g$ arbitrary, $x:=w_{k+1}$ and $y:=\bar{w}_k$ followed by Proposition \ref{prp_1}(\ref{prp_1_itm_1}), we get
    \allowdisplaybreaks
    \begin{align}
            &\lambda_k (\langle A(w_k), w_{k+1}-w\rangle + g(w_{k+1})-g(w))\nonumber\\&\leq\langle\nabla h(\bar{w}_k)-\nabla h(w_{k+1}), w_{k+1}-w\rangle\nonumber\\&
            =B_{h}(w,\bar{w}_k)-B_{h}(w,w_{k+1})-B_{h}(w_{k+1}, \bar{w}_k).\label{eqn_lmprf_2}
    \end{align}
    \allowdisplaybreaks
    By putting $k=k-1,w=w_{k+1}$ in (\ref{eqn_lmprf_2}) and using $\phi\nabla h(\bar{w}_k)=(\phi-1)\nabla h(w_k)+\nabla h(\bar{w}_{k-1})$ followed by Proposition \ref{prp_1}(\ref{prp_1_itm_1}), we get
    \allowdisplaybreaks
    \begin{align}
            &\lambda_{k-1}(\langle A(w_{k-1}), w_k-w_{k+1}\rangle + g(w_k)-g(w_{k+1}))\nonumber
            \\&\leq \langle \nabla h(\bar{w}_{k-1})-\nabla h(w_k), w_k-w_{k+1}\rangle\nonumber\\&=\phi\langle\nabla h(\bar{w}_k)-\nabla h(w_k), w_k-w_{k+1}\rangle\nonumber\\&=\phi\Big(B_{h}(w_{k+1},\bar{w}_k)-B_{h}(w_{k+1},w_k)-B_{h}(w_k,\bar{w}_k)\Big).\label{eqn_lmprf_3}
    \end{align}
On multiplying both sides of (\ref{eqn_lmprf_3}) by $\frac{\lambda_k}{\lambda_{k-1}}$, we get
\allowdisplaybreaks
\begin{align}
    &\lambda_k(\langle A(w_{k-1}),w_k-w_{k+1}\rangle +g(w_k)-g(w_{k+1}))\nonumber\\&\leq \frac{\lambda_k}{\lambda_{k-1}}\phi\Big(B_h(w_{k+1},\bar{w}_k)-B_h(w_{k+1},w_k)-B_h(w_k,\bar{w}_k)\Big).\label{eqn_mp_new_1}
\end{align}
By Lemma \ref{mm_lm}(\ref{lm_prt_3}), there exists $k_1\in\mathbb{N}$ such that for all $k\geq k_1$
\[\lambda_k=(1+\gamma_{k-1})\lambda_{k-1}.\]
Also, since $\displaystyle\sum_{k=0}^\infty \gamma_k < +\infty$, there exists $k_2\in\mathbb{N}$ such that
\[\gamma_{k-1}\leq \frac{1}{2\phi^2}~\forall k\geq k_2.\]
Therefore, by taking $k''=\max\{k_1,k_2\}$, we get
\[\frac{\lambda_k}{\lambda_{k-1}}\phi\leq \left(\phi+\frac{1}{2\phi}\right)~\forall k\geq k''.\]
Thus, we can rewrite (\ref{eqn_mp_new_1}) as
\allowdisplaybreaks
\begin{align}
    &\lambda_k(\langle A(w_{k-1}),w_k-w_{k+1}\rangle +g(w_k)-g(w_{k+1}))\nonumber\\&\leq \left(\phi + \frac{1}{2\phi}\right)\Big(B_h(w_{k+1},\bar{w}_k)-B_h(w_{k+1},w_k)-B_h(w_k,\bar{w}_k)\Big)~\forall k\geq k''.\label{eqn_mp_new_2}
\end{align}
 Let $w^{*}\in \mathcal{S}$. By setting $w=w^*$ in (\ref{eqn_lmprf_2}) and summing with (\ref{eqn_mp_new_2}), and then rearranging yields
 \allowdisplaybreaks
 \begin{align}
     &\lambda_k(\langle A(w_k),w_{k}-w^*\rangle +g(w_{k})-g(w^*))\nonumber\\ &\leq B_h(w^*,\bar{w}_k)-B_h(w^*,w_{k+1})-B_h(w_{k+1},\bar{w}_k)\nonumber\\&~~~~+\left(\phi + \frac{1}{2\phi}\right)\Big(B_h(w_{k+1},\bar{w}_k)-B_h(w_{k+1},w_k)-B_h(w_k,\bar{w}_k)\Big)\nonumber\\&~~~~+\lambda_k\langle A(w_k)-A(w_{k-1}),w_k-w_{k+1}\rangle~\forall k\geq k''\label{eqn_mp_new_3}.
 \end{align}
 By (\ref{eqn_defvi}) and the monotonicity of $A$, we get 
    \begin{equation}\label{eqn_lmprf_5}
        \begin{split}
            0&\leq \langle A(w^*),w_k-w^*\rangle+g(w_k)-g(w^*)\leq \langle A(w_k),w_k-w^*\rangle + g(w_k)-g(w^*). 
        \end{split}
    \end{equation}
    Combining (\ref{eqn_mp_new_3}) and (\ref{eqn_lmprf_5}) followed by Lemma \ref{lm_mm2}, we get
    \allowdisplaybreaks
    \begin{align}
        B_h(w^*,w_{k+1})&\leq B_h(w^*,\bar{w}_k)+\left(\phi+\frac{1}{2\phi}-1\right)B_h(w_{k+1},\bar{w}_k)\nonumber\\&~~~-\left(\phi+\frac{1}{2\phi}\right)B_h(w_k,\bar{w}_k)-\left(\frac{\phi}{2}+\frac{1}{2\phi}\right)B_h(w_{k+1},w_k)\nonumber\\&~~~+\frac{\phi}{2}B_h(w_k,w_{k-1})~\forall k\geq k':=\max\{\bar{k},k''\}.\label{eqn_mp_new_4}
    \end{align}
    From (\ref{eqn_stp2_1}), we have
    \[\nabla h(w_{k+1})=\frac{\phi}{\phi-1}\nabla h(\bar{w}_{k+1})-\frac{1}{\phi-1}\nabla h(\bar{w}_k)=(\phi+1)\nabla h(\bar{w}_{k+1})-\phi\nabla h(\bar{w}_k).\]
    Now, by applying Proposition \ref{prp_1}(\ref{prp_1_itm_2}), we get
    \begin{equation}\label{eqn_lmprf_8_new}
       \begin{split}
            (\phi +1)B_{h}(w^*,\bar{w}_{k+1})&=B_{h}(w^*,w_{k+1})+(\phi +1)B_{h}(w_{k+1},\bar{w}_{k+1})\\&~~~+\phi\Big(B_{h}(w^*,\bar{w}_k)-B_{h}(w_{k+1},\bar{w}_k)\Big).
       \end{split}
    \end{equation}
    On combining (\ref{eqn_mp_new_4}) and (\ref{eqn_lmprf_8_new}), the following is true for all $k\geq k'$
    \allowdisplaybreaks
    \begin{align}
        &(\phi +1)B_h(w^*,\bar{w}_{k+1})+\frac{\phi}{2}B_h(w_{k+1},w_k)-\phi B_h(w_{k+1},\bar{w}_{k+1})\nonumber\\&\leq (\phi +1)B_h(w^*,\bar{w}_k)+\frac{\phi}{2}B_h(w_k,w_{k-1})-\phi B_h(w_k,\bar{w}_k)+B_h(w_{k+1},\bar{w}_{k+1})\nonumber\\&~~~~+\left(\frac{1}{2\phi}-1\right)B_h(w_{k+1},\bar{w}_k).\label{eqn_mp_new_5}
    \end{align}
    From (\ref{eqn_stp2_1}), $\nabla h(\bar{w}_{k+1})=\frac{\phi-1}{\phi}\nabla h(w_{k+1})+\frac{1}{\phi}\nabla h(\bar{w}_k)$. Proposition \ref{prp_1}(\ref{prp_1_itm_2}) gives
\allowdisplaybreaks
\begin{align}
  B_{h}(w_{k+1},\bar{w}_{k+1})&=\frac{\phi-1}{\phi}\Big(B_{h}(w_{k+1},w_{k+1})-B_{h}(\bar{w}_{k+1},w_{k+1})\Big)\nonumber\\&~~~+\frac{1}{\phi}\Big(B_{h}(w_{k+1},\bar{w}_k)-B_{h}(\bar{w}_{k+1},\bar{w}_k)\Big)\nonumber\\&\leq\frac{1}{\phi} B_{h}(w_{k+1},\bar{w}_k). \label{eqn_lmprf_10} 
\end{align}
On combining (\ref{eqn_mp_new_5}) and (\ref{eqn_lmprf_10}), we get the second inequality in (\ref{eqn_lmprf_main}).\\
Next, we show the first inequality in (\ref{eqn_lmprf_main}). Applying Proposition \ref{prp_1}(\ref{prp_1_itm_1}) with $w=w^*, x=\bar{w}_{k+1}, y=w_{k+1}$, followed by (\ref{eqn_lmprf_2}), we get
\allowdisplaybreaks
\begin{align}
    &B_{h}(w^*,w_{k+1})+B_{h}(w_{k+1},\bar{w}_{k+1})\nonumber\\&=B_{h}(w^*,\bar{w}_{k+1})+\langle\nabla h(\bar{w}_{k+1})-\nabla h(w_{k+1}),w^*-w_{k+1}\rangle\nonumber\\&=B_{h}(w^*,\bar{w}_{k+1})+\frac{1}{\phi}\langle\nabla h(\bar{w}_k)-\nabla h(w_{k+1}),w^*-w_{k+1}\rangle\nonumber\\&\leq B_{h}(w^*,\bar{w}_{k+1})+\frac{\lambda_k}{\phi}\Big(\langle A(w_k),w^*-w_{k+1}\rangle-g(w_{k+1})+g(w^*)\Big)\nonumber\\&=B_{h}(w^*,\bar{w}_{k+1})+\frac{\lambda_k}{\phi}\langle A(w_k)-A(w_{k+1}),w^*-w_{k+1}\rangle\nonumber\\&~~~~+\frac{\lambda_k}{\phi}\Big(\langle A(w_{k+1}),w^*-w_{k+1}\rangle-g(w_{k+1})+g(w^*)\Big)\nonumber\\&\leq B_{h}(w^*,\bar{w}_{k+1})+\frac{\lambda_k}{\phi}\langle A(w_k)-A(w_{k+1}),w^*-w_{k+1}\rangle.\label{eqn_lmprf_11}
\end{align}
    By Cauchy-Schwarz inequality, Lemma \ref{mm_lm}(\ref{lm_prt_3}), and $\alpha$-strong convexity of $h$, it follows that for all $k\geq k_1$,
    \begin{equation}\label{eqn_lmprf_12}
        \begin{split}
            &\frac{\lambda_k}{\phi}\langle A(w_k)-A(w_{k+1}),w^*-w_{k+1}\rangle\\&\leq \frac{\lambda_k}{\phi}\lVert A(w_k)-A(w_{k+1})\rVert~\lVert w^*-w_{k+1}\rVert\\&\leq \frac{\lambda_{k}}{\phi}\frac{\eta_0\alpha}{\lambda_k}\lVert w_{k+1}-w_k\rVert~\lVert w^*-w_{k+1}\rVert\\&\leq\frac{1}{2}\frac{\alpha}{2}\Big(\lVert w_{k+1}-w_k\rVert^{2}+\lVert w^*-w_{k+1}\rVert^{2}\Big)\\&\leq \frac{1}{2}\Big(B_{h}(w_{k+1},w_k)+B_{h}(w^*,w_{k+1})\Big).
        \end{split}
    \end{equation}
    From (\ref{eqn_lmprf_11}) and (\ref{eqn_lmprf_12}), for all $k\geq k'$,  the following inequality holds
    \begin{equation}\label{eqn_roc_8}
        \begin{split}
            B_{h}(w_{k+1},\bar{w}_{k+1})&\leq B_{h}(w^*,\bar{w}_{k+1})+\frac{1}{2}B_{h}(w_{k+1},w_k)-\frac{1}{2}B_{h}(w^*,w_{k+1})\\&\leq B_{h}(w^*,\bar{w}_{k+1})+\frac{1}{2}B_{h}(w_{k+1},w_k),
        \end{split}
    \end{equation}
    and thus for all $k \geq k'$, we get
    \begin{equation*}
        \begin{split}
            0\leq&~B_{h}(w^*,\bar{w}_{k+1})+\frac{1}{2}B_{h}(w_{k+1},w_k)-B_{h}(w_{k+1},\bar{w}_{k+1})\\\implies 0 \leq&~ (\phi+1)B_{h}(w^*,\bar{w}_{k+1})+\frac{\phi}{2}B_{h}(w_{k+1},w_k)-\phi B_{h}(w_{k+1},\bar{w}_{k+1}),
        \end{split}
    \end{equation*}
    which is the first inequality of (\ref{eqn_lmprf_main}), and thus the proof is complete.  
\qed
\begin{theorem}
    Suppose that Assumptions A.\ref{assum_1}-A.\ref{assum_4} hold. Then, the sequences $(w_k)$ and $(\bar{w}_k)$ generated by Algorithm \ref{algo_1} converge to a point in $\mathcal{S}$.
\end{theorem}
\allowdisplaybreaks
{\it \textbf{Proof}}
    Let $w^*\in\mathcal{S}$ be arbitrary and let $(\mu_k)$ denote the sequence given by 
    \begin{equation}\label{eqn_roc_6}
        \mu_k:=(\phi +1)B_{h}(w^*,\bar{w}_k)+\frac{\phi}{2}B_{h}(w_k,w_{k-1})-\phi B_{h}(w_k,\bar{w}_k)~\forall k\in\mathbb{N}.
    \end{equation}
    From Lemma \ref{lm_1}, it is evident that there exists $k'\in\mathbb{N}$ such that the following is true for all $k\geq k'$
    \[0\leq \mu_{k+1}\leq\mu_k -\left(1-\frac{3}{2\phi}\right)B_h(w_{k+1},\bar{w}_k).\] Therefore, by Lemma \ref{lm_pre_1}, it follows that $\lim_{k\rightarrow\infty}\mu_k$ exists and $B_{h}(w_{k+1},\bar{w}_k)\rightarrow 0$ as $k\rightarrow\infty.$ Then, (\ref{eqn_lmprf_10}) implies that $B_{h}(w_{k+1},\bar{w}_{k+1})\rightarrow 0$. Also, by Proposition \ref{prp_1}(\ref{prp_1_itm_2}) with $\nabla h(w_k)=(\phi +1)\nabla h(\bar{w}_k)-\phi\nabla h(\bar{w}_{k-1})$, we get
    \allowdisplaybreaks
    \begin{align*}
            B_{h}(w_{k+1},w_k)&=(\phi +1)\Big(B_{h}(w_{k+1}, \bar{w}_k)-B_{h}(w_k,\bar{w}_k)\Big)-\phi \Big(B_{h}(w_{k+1},\bar{w}_{k-1})\\&~~~-B_{h}(w_k,\bar{w}_{k-1})\Big)\\&\leq (\phi +1)B_{h}(w_{k+1},\bar{w}_k)+\phi B_{h}(w_k,\bar{w}_{k-1})\rightarrow 0.
    \end{align*}
    Thus, we get
    \[\lim_{k\rightarrow\infty}\mu_k=(\phi +1)\lim_{k\rightarrow \infty}B_{h}(w^*,\bar{w}_k),\]
    and, in particular, $\lim_{k\rightarrow\infty}B_{h}(w^*,\bar{w}_k)$ exists.\\
    Using $\alpha$-strong convexity of $h$, we conclude that $w_{k+1}-\bar{w}_k\rightarrow 0$ and that $(w_k)$ and $(\bar{w}_k)$ are bounded. Thus, let $\bar{w}\in \mathcal{H}$ be a cluster point of $(\bar{w}_k)$. Then, there exists a subsequence $(\bar{w}_{{k}_{j}})$ such that $\bar{w}_{{k}_{j}}\rightarrow \bar{w}$ and $w_{k_{j}+1}\rightarrow\bar{w}$ as $j\rightarrow\infty$. Referring back to (\ref{eqn_lmprf_2}), for all $w\in \mathcal{H}$, we obtain
    \begin{multline*}
        \lambda_k(\langle A(w_{{k}_{j}}), w_{k_{j}+1}-w\rangle+ g(w_{k_{j}+1})-g(w))\\\leq \langle\nabla h(\bar{w}_{{k}_{j}})-\nabla h(w_{k_{j}+1}), w_{k_{j}+1}-w\rangle,
    \end{multline*}
    and taking the limit infimum on both sides as $j\rightarrow \infty$ shows that $\bar{w}\in\mathcal{S}$. Since $w^*\in\mathcal{S}$ in Lemma \ref{lm_1} was arbitrary, we now set $w^*=\bar{w}$. This implies that $\lim_{j\rightarrow\infty}B_{h}(w^*,\bar{w}_{{k}_{j}})=0$, leading to $\lim_{j\rightarrow\infty}\mu_{{k}_{j}}=0$. Also, from Lemma \ref{lm_1}, for $n\geq k_j$, we have $\mu_{n}\leq\mu_{{k}_{j}}$, and thus,
    \[(\phi+1)\lim_{n\rightarrow\infty}B_{h}(w^*,\bar{w}_n)=\lim_{n\rightarrow\infty}\mu_n\leq\lim_{j\rightarrow\infty}\mu_{{k}_{j}}=0,\]
    and thus $\bar{w}_k\rightarrow w^*$ from strong convexity of $h$. The convergence $w_k\rightarrow w^*$ follows from the fact that $w_k-\bar{w}_k\rightarrow 0.
    $
    \qed
\subsection{Linear Convergence of Modified B-GRAAL}
We next establish the R-linear convergence rate of Algorithm \ref{algo_1}. For this, we assume that for some $\sigma>0$,
\begin{equation}\label{eqn_roc_1}
    \langle A(z)-A(w),z-w\rangle\geq \sigma B_h(z,w)~\forall (z,w)\in\dom h\times \inte\dom h.
\end{equation}
Observe that when $h=\frac{1}{2}\lVert .\rVert^2$, (\ref{eqn_roc_1}) simplifies to the standard definition of $\sigma$-strong monotonicity. 
\allowdisplaybreaks
\begin{theorem}
    Assume that Assumptions A.\ref{assum_1}-A.\ref{assum_4} hold. Furthermore, assume that condition (\ref{eqn_roc_1}) is met. Then, the sequences $(w_k)$ and $(\bar{w}_k)$ generated by Algorithm \ref{algo_1} converge R-linearly to the unique point in $\mathcal{S}$.
\end{theorem}
\allowdisplaybreaks
{\it \textbf{Proof}}~
Without loss of generality, we assume that $\sigma<\frac{1}{2}$. Also, in the definition of $\lambda_k$ by multiplying right hand sides of (\ref{eqn_roc_2}) and (\ref{eqn_roc_3}) by $\sqrt{1-\sigma}$, and by replacing these updated values of $\lambda_k$ in the proof of Lemma \ref{lm_mm2}, we get the following updated bound.
\allowdisplaybreaks
\begin{align}\label{eqn_roc_4}
   &\lambda_k\langle A(w_k)-A(w_{k-1}), w_k-w_{k+1}\rangle\nonumber\\&\leq \frac{\phi}{2}(1-\sigma)B_h(w_k,w_{k-1})+\frac{\phi}{2}B_h(w_{k+1},w_k)~\forall k\geq \bar{k}.
\end{align}
Let $w^*\in\mathcal{S}$. Since (\ref{eqn_roc_1}) holds, we get
\allowdisplaybreaks
\begin{align}\label{eqnn}
        0&\leq \langle A(w^*), w_k-w^*\rangle + g(w_k)-g(w^*)\nonumber\\&\leq \langle A(w_k),w_k-w^*\rangle+g(w_k)-g(w^*)-\sigma B_h(w^*,w_k).
\end{align}
Following a similar approach to the proof of Lemma \ref{lm_1}, with (\ref{eqnn}) replacing (\ref{eqn_lmprf_5}) and (\ref{eqn_roc_4}) replacing Lemma \ref{lm_mm2}, the corresponding form of (\ref{eqn_mp_new_4}) becomes
\allowdisplaybreaks
\begin{align}\label{eqn_roc_7}
         &B_h(w^*,w_{k+1})\leq B_h(w^*,\bar{w}_k)+(\phi-1)B_h(w_{k+1},\bar{w}_k)-\phi B_h(w_k,\bar{w}_k)\nonumber\\&-\sigma B_h(w^*,w_k)
    -\frac{\phi}{2}B_h(w_{k+1},w_k)+\frac{(1-\sigma)\phi}{2}B_h(w_k,w_{k-1})~\forall k\geq k'.
\end{align}
Proposition \ref{prp_1}(\ref{prp_1_itm_2}) with $\nabla h(w_{k+1})=(\phi+1)\nabla h(\bar{w}_{k+1})-\phi\nabla h(\bar{w}_k)$ gives
\begin{equation}\label{eqn_roc_5}
    \begin{split}
        B_h(w^*,w_{k+1})=&(\phi+1)\Big(B_h(w^*,\bar{w}_{k+1})-B_h(w_{k+1},\bar{w}_{k+1})\Big)\\&-\phi\Big(B_h(w^*,\bar{w}_k)-B_h(w_{k+1},\bar{w}_k)\Big).
    \end{split}
\end{equation}
Let $\mu_k$ represent the sequence defined by (\ref{eqn_roc_6}). Then, applying (\ref{eqn_roc_5}) to substitute $B_h(w^*,w_{k+1})$ and $B_h(w^*,w_{k})$ in (\ref{eqn_roc_7}), and subsequently using (\ref{eqn_lmprf_10}), we get the following for all $k\geq k'$,
\begin{align*}
    \mu'_{k+1}:=&~(\phi+1)B_h(w^*,\bar{w}_{k+1})+\frac{\phi}{2}B_h(w_{k+1},w_k)-B_h(w_{k+1},\bar{w}_{k+1})\\\leq&~ (1-\sigma)(\phi+1)B_h(w^*,\bar{w}_k)+\frac{(1-\sigma)\phi}{2}B_h(w_k,w_{k-1})+\\&~(\sigma(1+\phi)-\phi)B_h(w_k,\bar{w}_k)+\sigma\phi B_h(w^*,\bar{w}_{k-1})-\sigma\phi B_h(w_k,\bar{w}_{k-1})\\&~+\phi B_h(w_{k+1},\bar{w}_{k+1})+\left(\frac{1}{2\phi}-1\right)B_h(w_{k+1},\bar{w}_k)\\\leq&~ (1-\sigma)\mu_k+\sigma\phi B_h(w^*,\bar{w}_{k-1}).
\end{align*}
Since $(\mu_k)$ is non increasing by Lemma \ref{lm_1} and $\frac{(1-\sigma)\phi}{1-\sigma(2-\phi)}\geq 1$, we get
\allowdisplaybreaks
\begin{align*}
        \mu'_{k+1}\leq&~(1-\sigma)\mu_{k-1}+\sigma(\phi-1)(\phi+1)B_h(w^*,\bar{w}_{k-1})\\=&~(1-\sigma(2-\phi))(\phi+1)B_h(w^*,\bar{w}_{k-1})+(1-\sigma)\frac{\phi}{2}B_h(w_{k-1},w_{k-2})\\&-(1-\sigma)\phi B_h(w_{k-1},\bar{w}_{k-1})\\\leq&~ (1-\sigma(2-\phi))\Big((\phi +1)B_h(w^*,\bar{w}_{k-1})+\frac{\phi}{2}B_h(w_{k-1},w_{k-2})\\&~-B_h(w_{k-1},\bar{w}_{k-1})\Big)\\=&~q\mu'_{k-1},~\text{where}~q=(1-\sigma(2-\phi))\in(0,1).
\end{align*}
Therefore, $\mu'_k$ converges Q-linearly to $0$. By (\ref{eqn_roc_8}) and $\alpha$-strong convexity of $h$, it follows that
\begin{align*}
        \mu'_{k+1}=&~\phi B_h(w^*,\bar{w}_{k+1})+\frac{\phi-1}{2}B_h(w_{k+1},w_k)\\&+\left(B_h(w^*,\bar{w}_{k+1})+\frac{1}{2}B_h(w_{k+1},w_k)-B_h(w_{k+1},\bar{w}_{k+1})\right)\\\geq&~\frac{\alpha\phi}{2}\lVert w^*-\bar{w}_{k+1}\rVert^2+\frac{\alpha(\phi-1)}{4}\lVert w_{k+1}-w_k\rVert^2.
\end{align*}
This shows that $\bar{w}_k$ converges to $w^*$ with R-linear convergence, and the difference $\lVert w_{k+1} - w_k \rVert$ also converges to zero at the same rate. As a result, the sequence $(w_k)$ exhibits R-linear convergence to $w^*$. Since the choice of $w^*$ from the set $\mathcal{S}$ was arbitrary, $w^*$ is unique.
\qed
\section{Numerical Results}\label{sec:5}
In this section, we compare our Algorithm \ref{algo_1} with B-GRAAL \citep[Algorithm 2]{tam2023bregman}, and B-aGRAAL \citep[Algorithm 3]{tam2023bregman}. Specifically, in Section \ref{sec5:2}, we present the comparison results for matrix game problem. In Section \ref{sec5:1}, we compare the algorithms in the context of the the sparse logistic regression problem. We run all the experiments in Python 3 on a Macbook Air with M1 chip and 8GB RAM. \par
Both the considered problems can be formulated as MVIP (\ref{eqn_defvi}) for different choices of the operator $A$ and the function $g$. Note that the solution of MVIP (\ref{eqn_defvi}) can also be formulated as the monotone inclusion: find $w^*\in\mathcal{H}$ such that
\[0\in A(w^*)+\partial g(w^*).\]
Thus, by (\ref{eqn_step_2}), we get
\[0\in A(w_k)+\partial g(w_{k+1})+\frac{1}{\lambda_k}(\nabla h(w_{k+1})-\nabla h(\bar{w}_k)),\]
and we track the quantity $J_k$ defined as
\begin{multline}
    J_k:=\frac{1}{\lambda_k}(\nabla h(\bar{w}_k)-\nabla h(w_{k+1}))+A(w_{k+1})-A(w_k)\in A(w_{k+1})+\partial g(w_{k+1}).
\end{multline}
as a natural residual for Algorithm \ref{algo_1}.
The parameters used for the execution of each algorithm are as follows.
\begin{enumerate}
\renewcommand{\labelenumi}{}
\item B-GRAAL: $\phi= \frac{\sqrt{5}+1}{2}, \lambda= \frac{\phi}{2L}.$
    \item B-aGRAAL: $\varphi=1.5$, $\lambda_0=\frac{\varphi}{2}\frac{\lVert w_1-w_0 \rVert}{\lVert A(w_1)-A(w_0)\rVert}$, $\rho= \frac{1}{\phi}+\frac{1}{{\phi}^2}$, $\lambda_{\max}=10^6.$ 
    \item Algorithm \ref{algo_1}: $\eta_1=0.75, \eta_0=0.80, \lambda_0=\frac{\phi}{2}\frac{\lVert w_1-w_0 \rVert}{\lVert A(w_1)-A(w_0)\rVert}$, and the sequence of positive terms $\gamma_{k}=\frac{r(\log (k+1))^s}{(k+1)^t}, r, s>0; t>1~\forall k\geq 1.$ 
\end{enumerate}
\subsection{Matrix Games}\label{sec5:2}
In this instance, we examine the following matrix game: 
\begin{equation}\label{eqn_app3_1}
    \min_{x\in\triangle ^k}\max_{y\in\triangle^k}\langle Px, y\rangle,
\end{equation}
where $\triangle^k:= \left\{x\in {\mathbb{R}}^{k}_+: \sum_{i=1}^k x_i=1\right\}$ denotes the unit simplex and $P\in {\mathbb{R}}^{k\times k}$ is a given matrix. Note that (\ref{eqn_app3_1}) can be formulated as a MVIP (\ref{eqn_defvi}).
Specifically, we examine the problem in the form (\ref{eqn_app3_1}), which involves positioning a server within a network $G = (U, E)$ consisting of  $k$ vertices, with the objective of minimizing its response time.
In this scenario, a request arises from an unknown vertex  $u_j \in U$, and the goal is to position the server at a vertex $u_i \in U$ in a manner that minimizes the response time, which is determined by the graphical distance $d(u_i, u_j)$.
We analyze the scenario where the request location  $u_j \in U$ is strategically determined by an adversary. The decision variable $x \in \triangle^k$ (similarly, $y \in \triangle^k$) represents mixed strategies for server placement (and request origin, respectively). More precisely, $x_t$ (respectively, $y_t$) denotes the probability of the server (or request origin) being positioned at node $u_t$ for $t = 1, 2, \ldots, k$. The matrix $P$ is defined as the distance matrix of the graph $G$, where each entry is given by $P_{ij} = d(u_i, u_j)$ for all vertices $( u_i, u_j) \in U$. Consequently, the objective function in (\ref{eqn_app3_1}) quantifies the expected response time, which we aim to minimize, while the adversary attempts to maximize it.\par
\begin{figure}[h]%
\centering
\mbox{\subfigure{\includegraphics[scale=0.38]{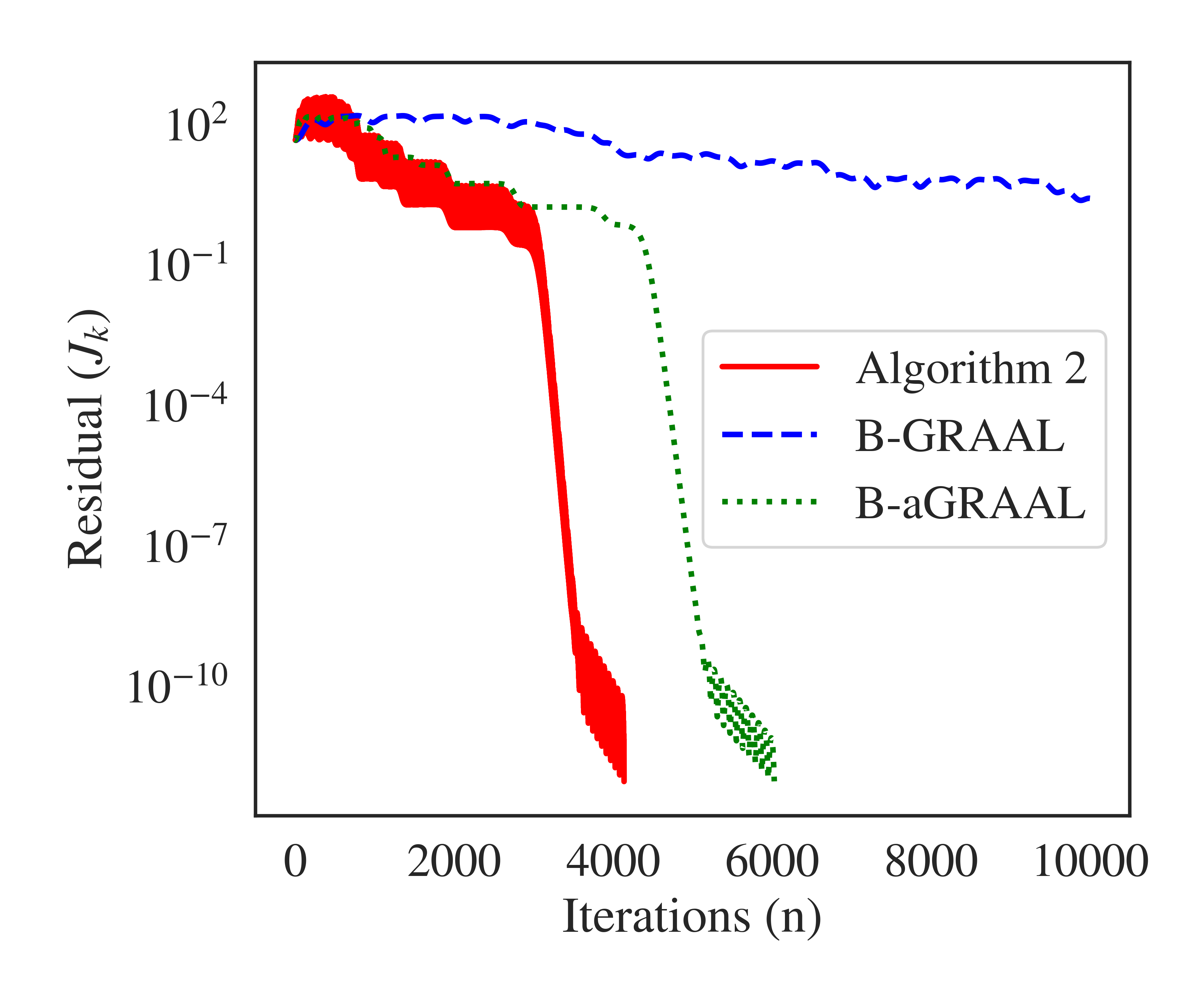} }\quad
\subfigure{{\includegraphics[scale=0.38]{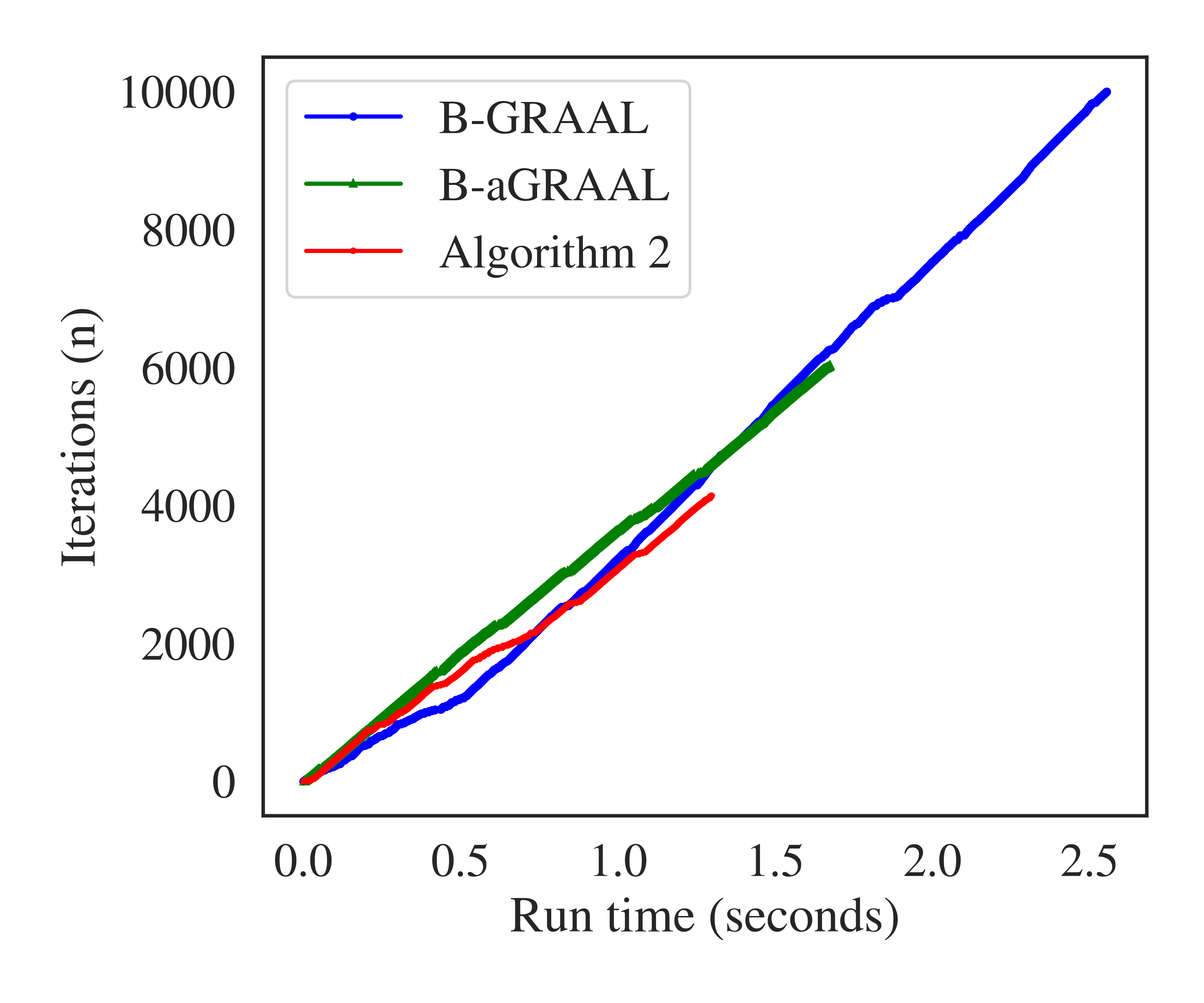} } }}
\caption{Matrix game results for $k=10$}
 \label{fig:n=10_matrixgame}
\end{figure}
\begin{figure}[h]
\centering
\mbox{\subfigure{\includegraphics[scale=0.38]{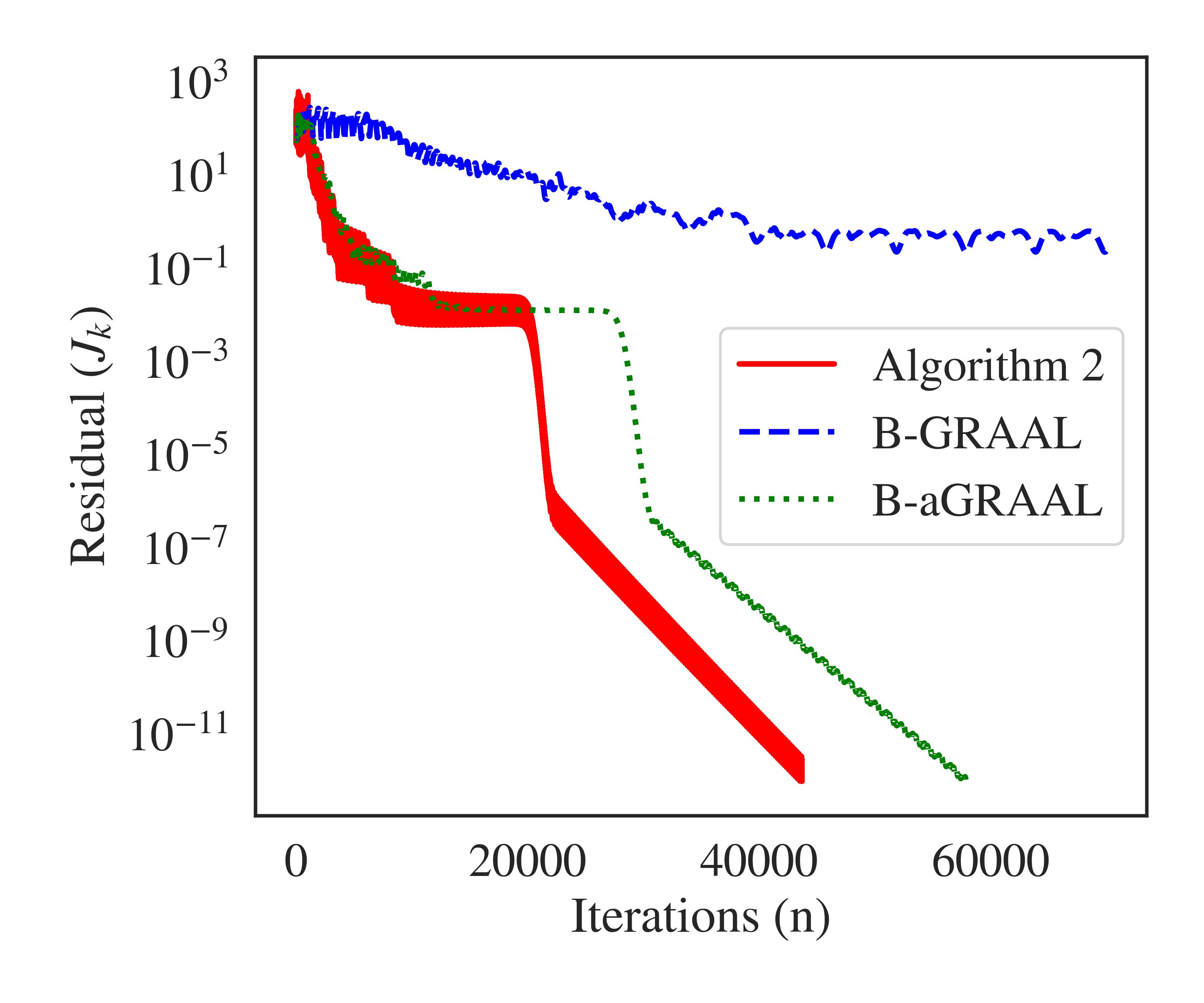} }\quad
\subfigure{{\includegraphics[scale=0.38]{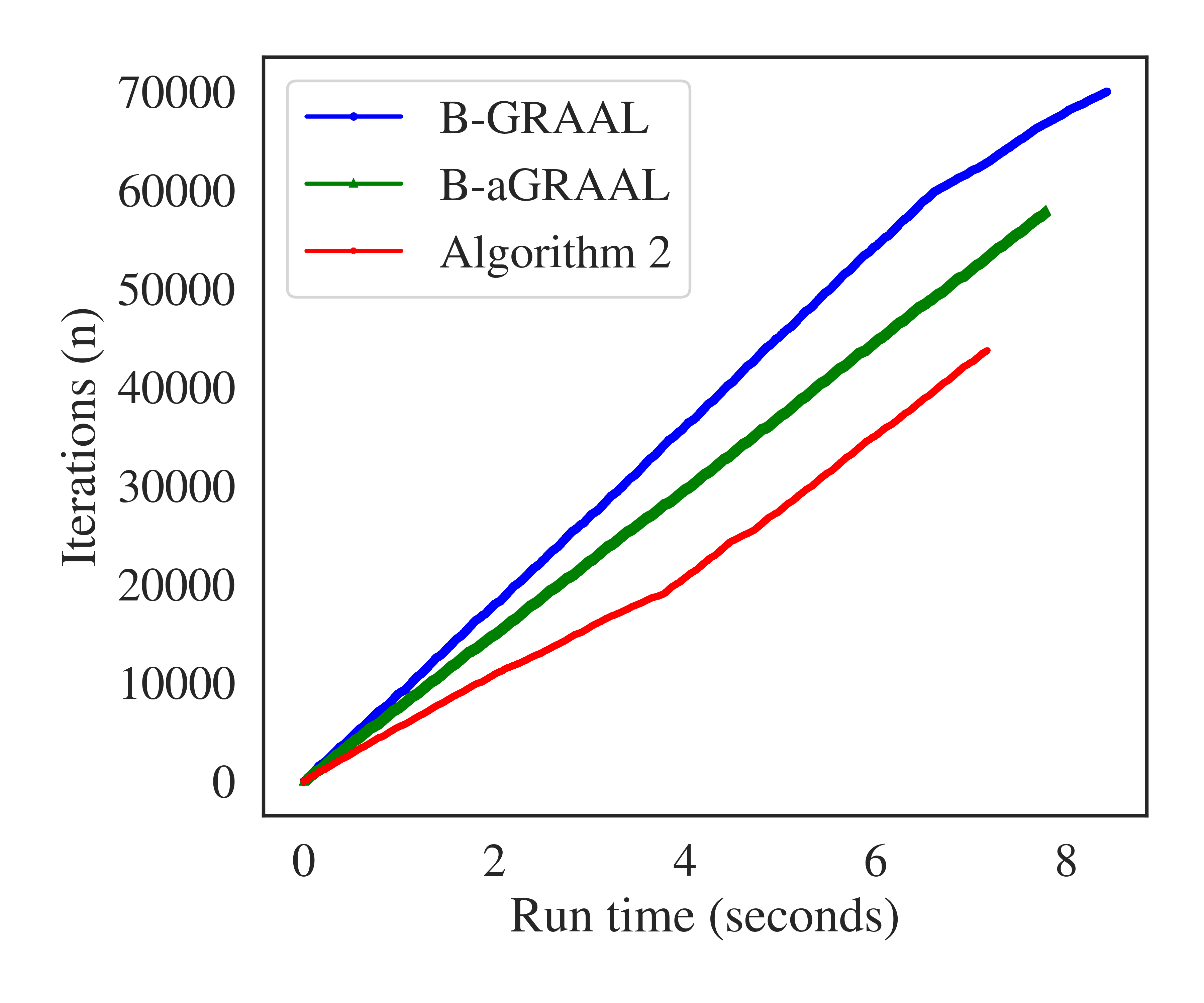} } }}
\caption{Matrix game results for $k=20$}
 \label{fig:n=20_matrixgame}
\end{figure}
\begin{figure}[h]%
\centering
\mbox{\subfigure{\includegraphics[scale=0.38]{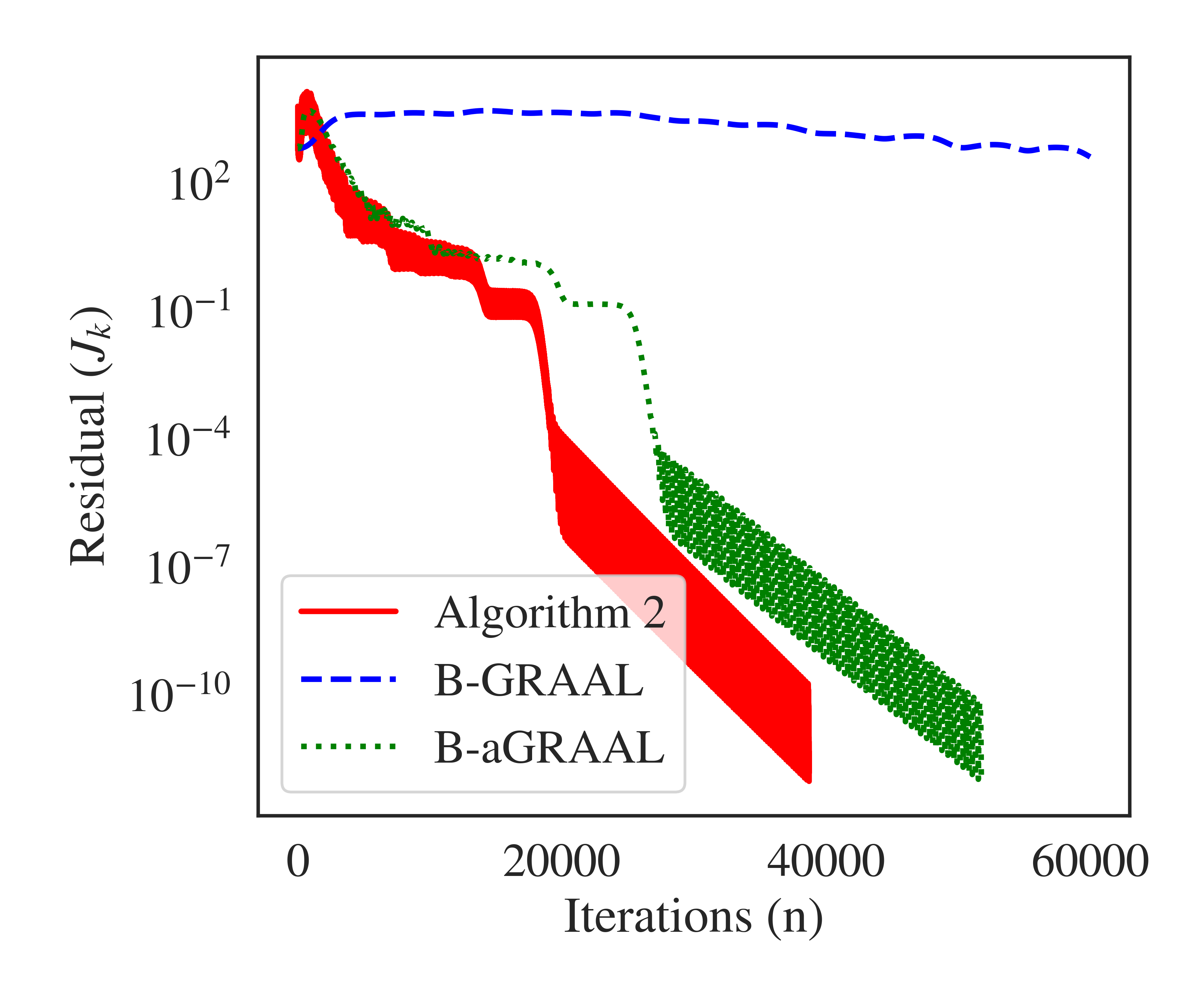} }\quad
\subfigure{{\includegraphics[scale=0.38]{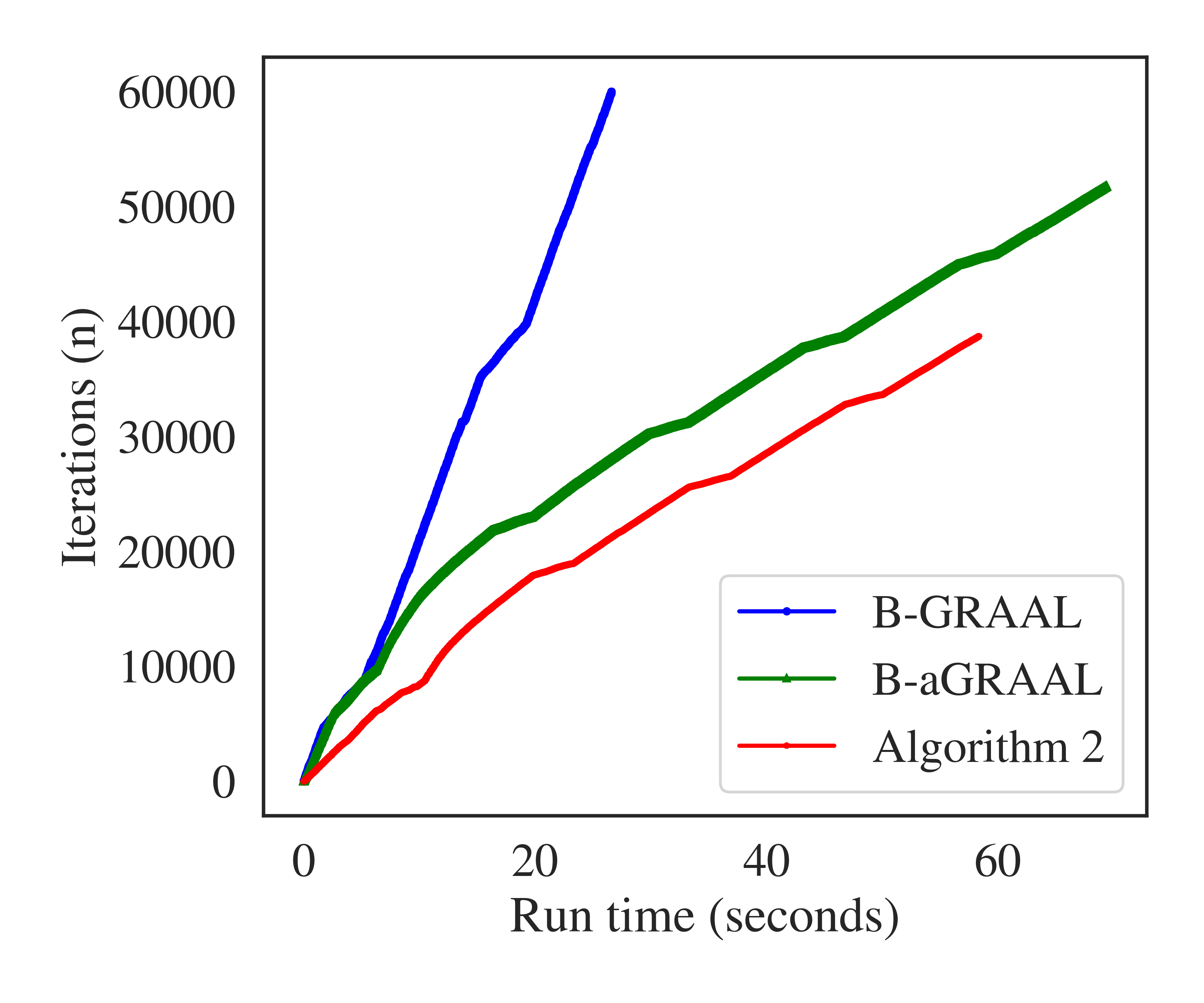} } }}
\caption{Matrix game results for $k=100$}
 \label{fig:n=100_matrixgame}
\end{figure}

We conducted three experiments with $k = 10$, $k= 20$ and $k = 100$, with the corresponding results presented in Figs. \ref{fig:n=10_matrixgame}, \ref{fig:n=20_matrixgame}, and \ref{fig:n=100_matrixgame}. The initial points are selected as $x_0 = y_0 = \left(\frac{1}{k}, \frac{1}{k}, \dots, \frac{1}{k}\right) \in \triangle^k$, with $w_0 = (x_0, y_0)$. Subsequently, $w_1$ is set as a randomly perturbed version of $w_0$. Lipschitz constant of $A$ is calculated as $L=\lVert P \rVert_2$ and $r=0.0007, s=7.5, t=1.1$. 
        \subsection{Sparse Logistic Regression}\label{sec5:1}
\begin{figure}[h]
\centering
\mbox{\subfigure{\includegraphics[scale=0.38]{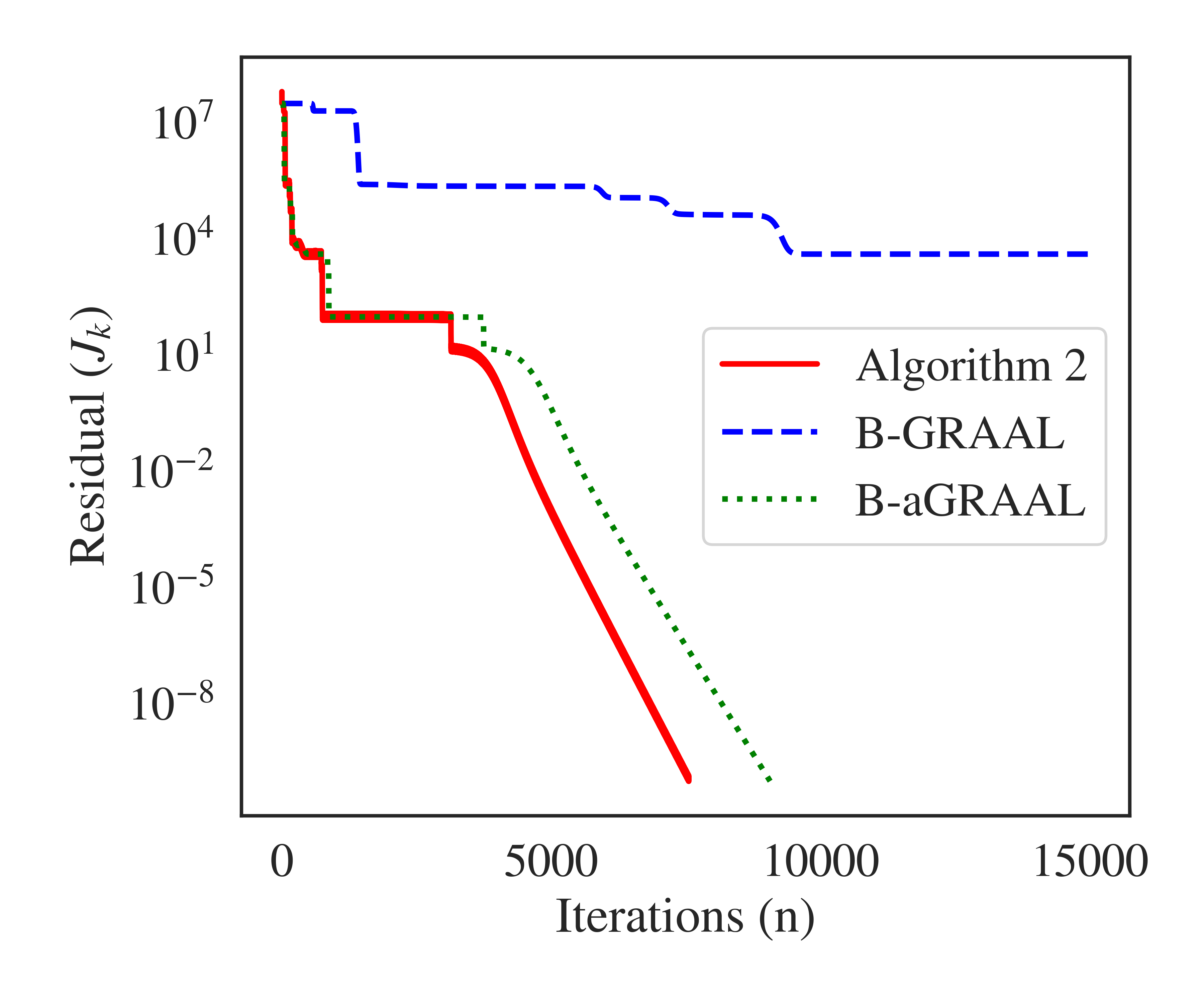} }\quad
\subfigure{{\includegraphics[scale=0.38]{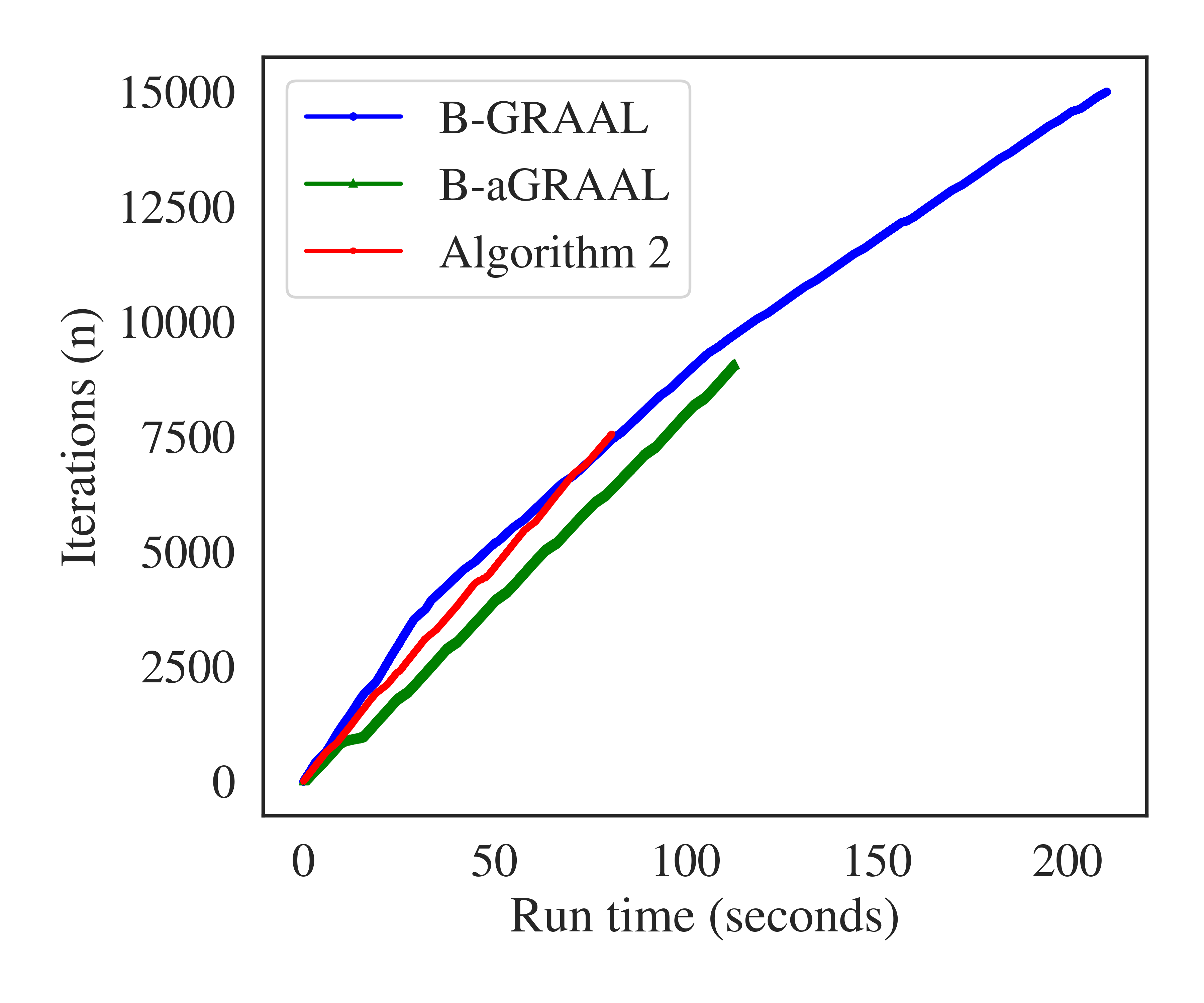} } }}
\caption{Result of the data ijcnn1.bz2 ($M=35000, n=22$)}
 \label{fig:n=slr_ijcnn1}
\end{figure}

\begin{figure}[h]
\centering
\mbox{\subfigure{\includegraphics[scale=0.38]{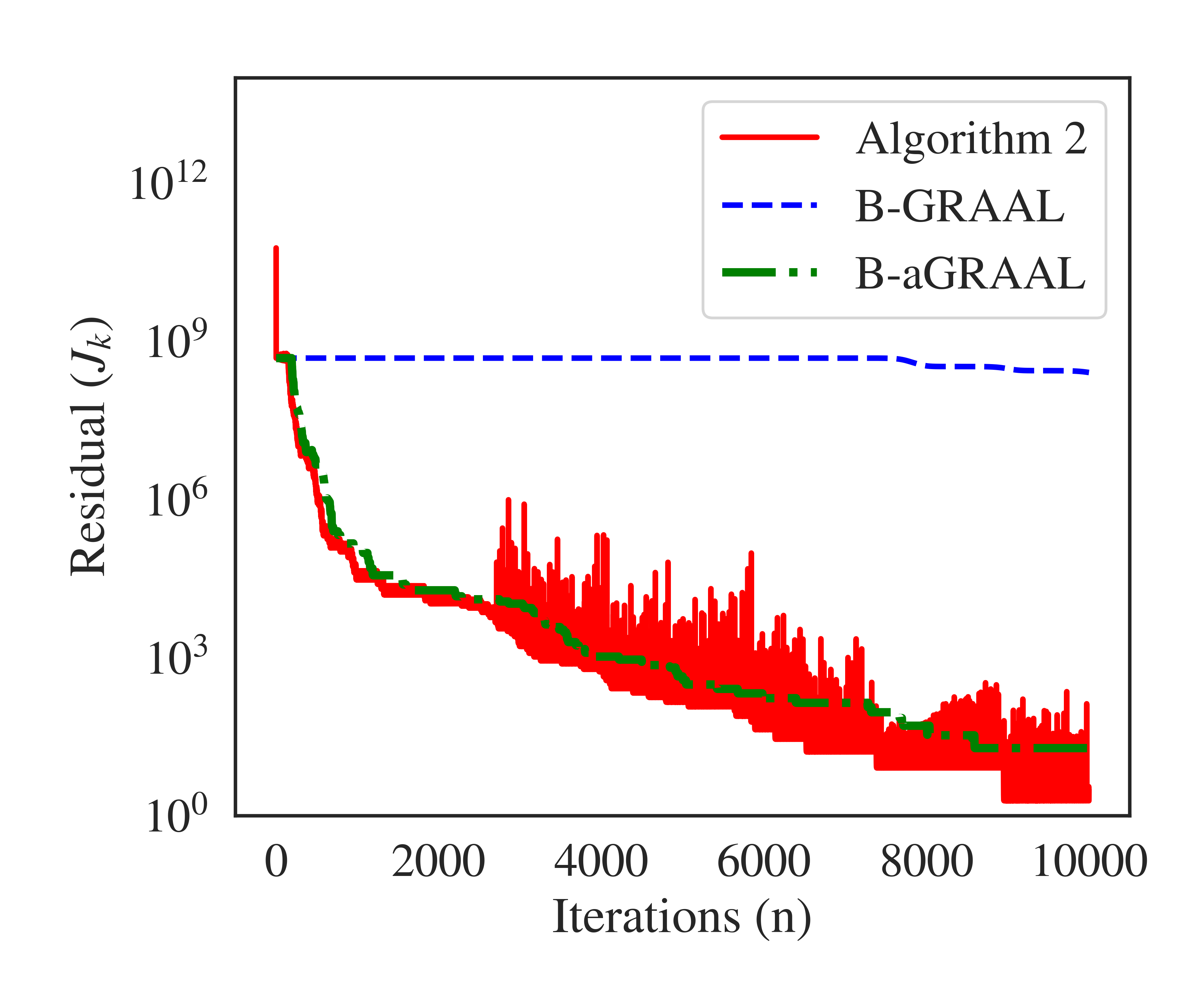} }\quad
\subfigure{{\includegraphics[scale=0.38]{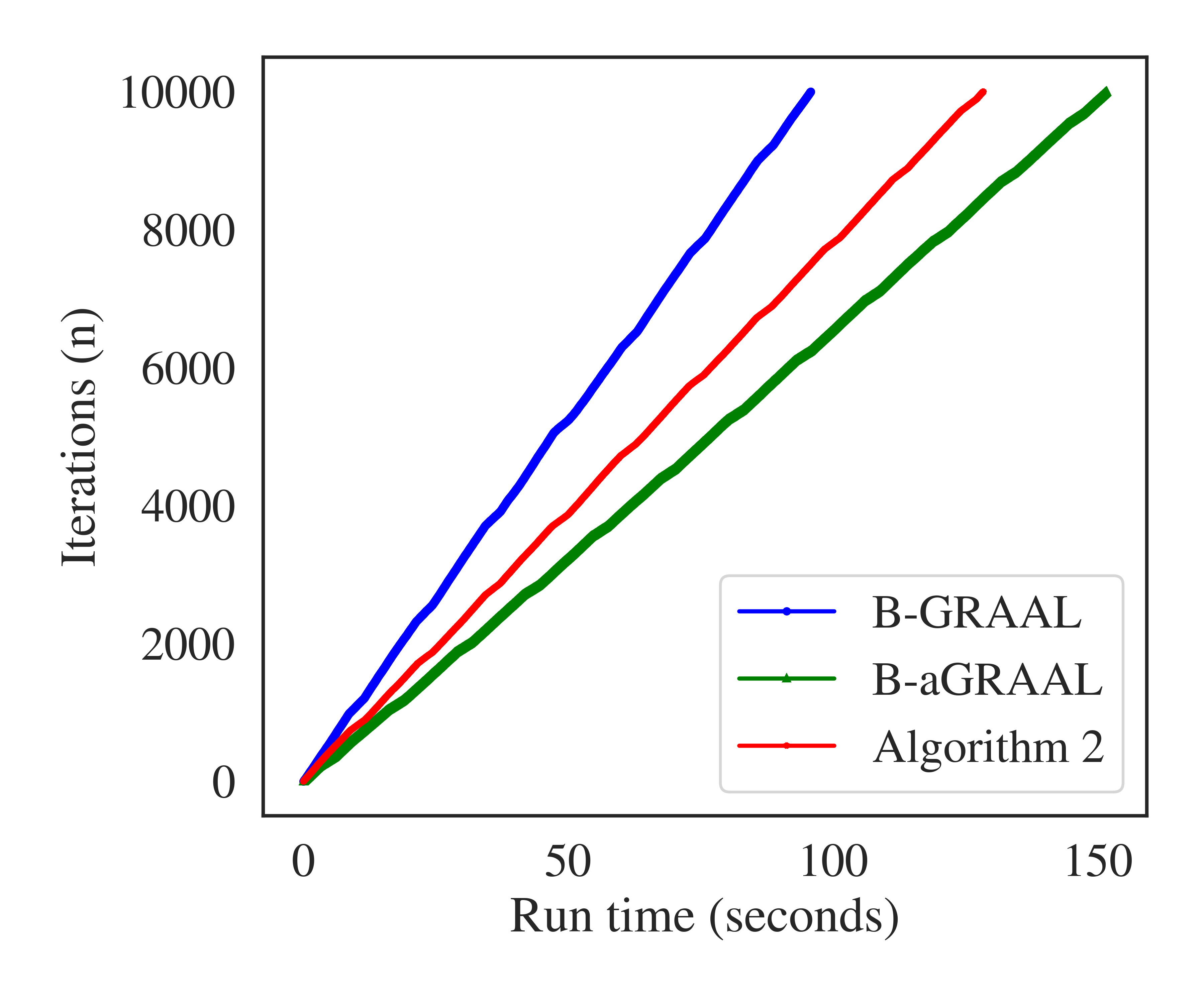} } }}
\caption{Result of the data a9a ($M=32561, n=123$)}
 \label{fig:n=slr_a9a}
\end{figure}

\begin{figure}[h]
\centering
\mbox{\subfigure{\includegraphics[scale=0.38]{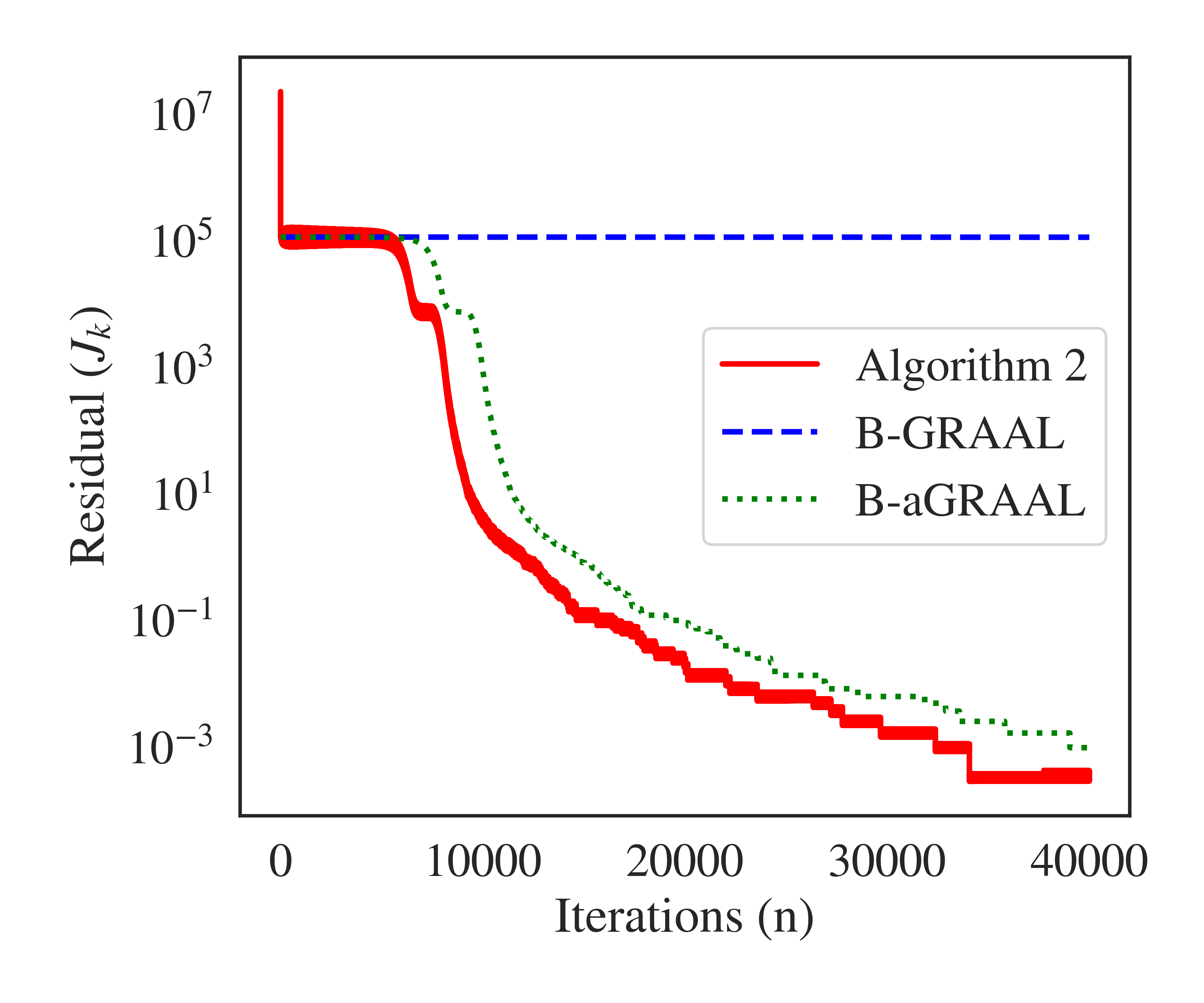} }\quad
\subfigure{{\includegraphics[scale=0.38]{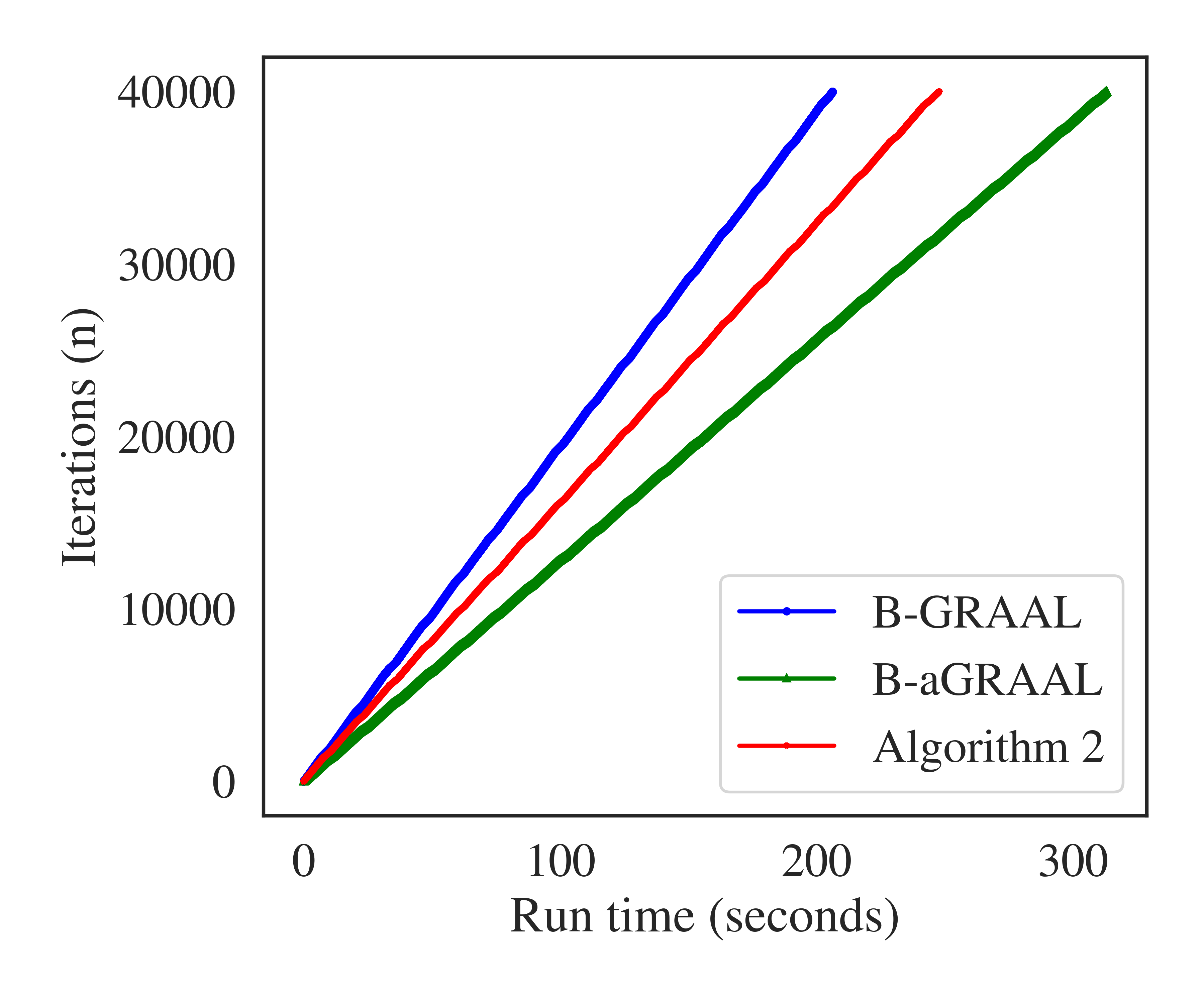} } }}
\caption{Result of the data duke.bz2 ($M=44, n=7129$)}
 \label{fig:n=slr_duke}
\end{figure}

In this problem, we are interested to find $x\in\mathbb{R}^n$ such that the following loss function is minimized
\[k(x)=\displaystyle\sum_{i=1}^M \log(1+\exp(-c_i\langle d_i,x\rangle))+\bar{\beta}\lVert x\rVert_1,\]
where $d^i\in{\mathbb{R}}^n, c_i\in\{-1,1\}, i=1,2,\ldots, M$ are observations; and $\bar{\beta}>0$ is a regularization parameter. It can be seen that (for details, see  \citep[Section 4]{hoai2024new}) the problem $\min_{x\in{\mathbb{R}}^n} k(x)$ can be formulated as problem (\ref{eqn_defvi}) with 
\[A(x)=\nabla h(x);~ h(x)=\displaystyle\sum_{i=1}^M \log(1+\exp(-c_i\langle d_i,x\rangle)); ~g(x)=\bar{\beta}\lVert x\rVert_1.\]
The tested datasets are sourced from the LIBSVM\footnote{https://www.csie.ntu.edu.tw/~cjlin/libsvmtools/datasets/} library, including ijcnn1.bz2, a9a, and duke.bz2. We choose $w_0=zeros(n), w_1=w_0+10^{-9}random(n),\\ \bar{\beta}=0.005\lVert C^{T}c\rVert_{\infty}(C~\text{composed by}~d^i, i=1,2,\ldots,M, c=(c_i)_{i=1}^N)$ for all algorithms. For the datasets ijcnn1 and a9a, we take $r=0.0001, s=7.2,$ and $ t=1.01$. For duke.bz2, $r=0.0005, s=7$ and $t=1.1$. The numerical results for the datasets ijcnn1.bz2, a9a, and duke.bz2 are presented in Figs. \ref{fig:n=slr_ijcnn1}, \ref{fig:n=slr_a9a}, and \ref{fig:n=slr_duke}, respectively.
\section{Conclusion}\label{sec:6}
In this article, we explored a modification of the B-GRAAL algorithm to solve mixed variational inequality problems. Under standard assumptions, we established the R-linear convergence of the modified algorithm. Unlike the B-GRAAL algorithm, the implementation of the modified algorithm does not require knowledge of the Lipschitz constant $L$ and the numerical experiments on the sparse logistic regression problem and the matrix game problem demonstrate that this modification significantly enhances the performance of the B-GRAAL algorithm.
\section*{Acknowledgements}
The authors wishes to thank Dr. P.T. Hoai and Dr. D. Uteda for sharing Python codes of their papers  \cite{hoai2024new} and \cite{tam2023bregman}, respectively. This helped the authors to write the Python codes for this paper.

\end{document}